\documentclass[11pt]{article}
\usepackage[T1]{fontenc}
\usepackage{lmodern}
\usepackage{amsmath,amssymb}
\usepackage{xcolor}
\usepackage{tikz}
\usepackage{tcolorbox}
\usepackage{geometry}
\usepackage{enumitem}
\usepackage{url}
\usepackage{hyperref}
\usetikzlibrary{
arrows.meta,
positioning,
shapes.geometric
}

\definecolor{mathblue}{RGB}{20,45,90}
\definecolor{journeygrayblue}{RGB}{70,85,120}
\definecolor{curiositycolor}{RGB}{75,105,160}
\definecolor{ancientmath}{RGB}{145,75,55}
\definecolor{origincolor}{RGB}{166,78,61}
\definecolor{countinggreen}{RGB}{55,120,85}
\definecolor{arithmeticorange}{RGB}{190,105,45}
\definecolor{numeralplum}{RGB}{105,65,125}
\definecolor{indianmathgold}{RGB}{175,120,35}
\definecolor{symbolicindigo}{RGB}{55,65,140}
\definecolor{sectionblue}{RGB}{0,70,140}
\definecolor{turquoise}{RGB}{40,180,170}
\definecolor{magenta}{RGB}{200,40,150}
\definecolor{rose}{RGB}{220,70,120}
\definecolor{storyblue}{RGB}{40,100,200}
\definecolor{forestgreen}{RGB}{40,140,80}
\definecolor{storyred}{RGB}{200,60,60}
\definecolor{mysterypurple}{RGB}{120,70,170}
\definecolor{magicgold}{RGB}{220,160,20}
\definecolor{oceanblue}{RGB}{30,120,180}
\definecolor{emerald}{RGB}{0,155,110}
\definecolor{crimson}{RGB}{180,30,60}
\definecolor{amber}{RGB}{240,170,40}
\definecolor{tealgreen}{RGB}{0,140,140}
\definecolor{indigo}{RGB}{75,0,130}
\definecolor{charcoal}{RGB}{60,60,60}
\title{
\textcolor{mathblue}{\Huge From Counting Objects to Negative Numbers}\\[4mm]
\textcolor{journeygrayblue}{\Large
The Journey to Understanding Why
\[(-)\times(-)=+\]}}

\author{
Chandradew Sharma\thanks{ \texttt{csharma@goa.bits-pilani.ac.in}}\\
Department of Physics\\
Birla Institute of Technology and Science, K K Birla Goa Campus\\
Zuarinagar, Goa 403726, India
}

\date{}

\begin{document}

\maketitle


\begin{abstract}
Why does $(-)\times(-)=+$? Although every student learns this rule, many accept it as a convention rather than as a mathematical necessity. This article develops a conceptual explanation by tracing the evolution of number from counting objects to representing opposites. Beginning with the emergence of numerals, place value, zero, and negative numbers, it shows how each stage extended the meaning of number and the operations defined on it. The discussion argues that multiplication of signed numbers is most naturally understood in terms of preserving or reversing direction, so that the familiar rule $(-)\times(-)=+$ follows as a logical consequence of a coherent arithmetic rather than as an arbitrary convention. By combining historical development with mathematical reasoning, the article offers an accessible framework for teachers and students seeking a deeper understanding of signed multiplication and illustrates how the evolution of mathematical ideas can enrich classroom learning.
\end{abstract}


\section*{\textcolor{curiositycolor}{The Question That Lingers}}

The teacher walked into the classroom, picked up a piece of chalk, and wrote a single expression on the board:
\[(-2)\times(-3)=\;?\]
Nobody spoke.
After a few moments, John raised his hand.
\begin{quote}
\emph{"I think the answer should be $-6$. If both numbers are negative, shouldn't the result stay negative?"}
\end{quote}
Hamid looked unconvinced.
\begin{quote}
\emph{"The textbook says the answer is $+6$. I know the rule because I memorized it, but I don't really understand why it works."}
\end{quote}
Rajesh kept looking at the expression.
\begin{quote}
\emph{"I understand
\[
2\times3,
\]
because it means two groups of three:
\[3+3=6.\]
But what does
\[(-2)\times(-3)\]
actually represent? What is a negative group of negative quantities? I know the rule, but I cannot picture the operation."}
\end{quote}
A smile appeared on the teacher's face.
\begin{quote}
\emph{"Those are exactly the questions that puzzled people for centuries."}
\end{quote}
Then the teacher continued.
\begin{quote}
\emph{"To understand why
\[(-)\times(-)=+,\]
we first have to ask a much older question: Why did numbers exist in the first place? Mathematics did not begin with symbols, formulas, or multiplication tables. It began with people trying to solve problems. Imagine a world without numerals, written arithmetic, or even a formal number system. People still needed answers. Which collection has more? Which has less? Are two collections the same size?"}
\end{quote}


\begin{figure}[htbp]
\centering

\begin{tikzpicture}[
node distance=0.7cm,
every node/.style={
draw,
rounded corners,
minimum width=2.7cm,
minimum height=0.95cm,
font=\small,
align=center
},
>=Stealth
]

\node[fill=yellow!25] (need)
{Human\\Need};

\node[fill=green!20,right=of need] (count)
{Counting};

\node[fill=cyan!20,right=of count] (arith)
{Arithmetic};

\node[fill=orange!25,right=of arith] (notation)
{Numbers};

\node[fill=violet!20,right=of notation] (place)
{Place\\Value};

\node[fill=blue!18,below=of place] (zero)
{Zero};

\node[fill=red!18,left=of zero] (negative)
{Negative\\Numbers};

\node[fill=teal!20,left=of negative] (signedarith)
{Extended \\Arithmetic};

\node[fill=magenta!20,left=of signedarith] (final)
{$(-)\times(-)=+$};

\draw[->,very thick] (need)--(count);
\draw[->,very thick] (count)--(arith);
\draw[->,very thick] (arith)--(notation);
\draw[->,very thick] (notation)--(place);
\draw[->,very thick] (place)--(zero);
\draw[->,very thick] (zero)--(negative);
\draw[->,very thick] (negative)--(signedarith);
\draw[->,very thick] (signedarith)--(final);

\end{tikzpicture}

\caption{
Conceptual progression showing how each new idea extends the limits of the previous one, culminating in \((-) \times (-)=+\).
}

\label{fig:roadmap}

\end{figure}
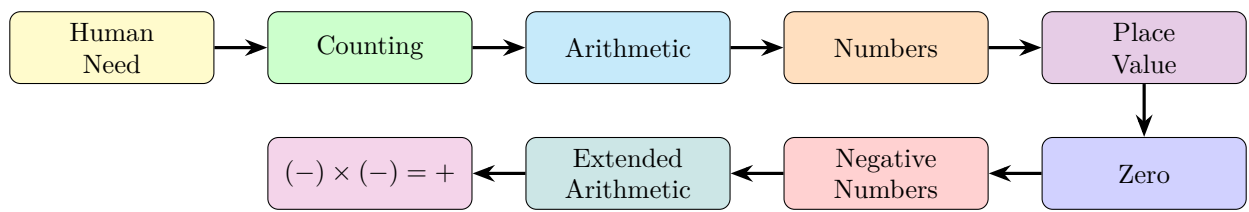

The teacher turned toward the diagram (Figure~\ref{fig:roadmap}) projected beside the board.
\begin{quote}
\emph{"This is the path we are going to follow . Every new idea grows out of a problem the previous one cannot solve. Counting leads to arithmetic. Arithmetic drives the need for written number systems.  Writing numbers gives rise to place value. Place value brings us to zero. Once zero is accepted, negative numbers are no longer far behind. Extending arithmetic to include those numbers eventually brings us back to the question we started with:
\[(-)\times(-)=+.\]
By the time we reach that point, the rule will not feel like something to memorize. It will simply be the answer that is consistent with everything that came before."}
\end{quote}

\section*{\textcolor{ancientmath}{The First Question: Which Collection is Larger?}}

The teacher looked around the classroom.
\begin{quote}
\emph{"Imagine a time before numerals, formal arithmetic, or written language."}
\end{quote}


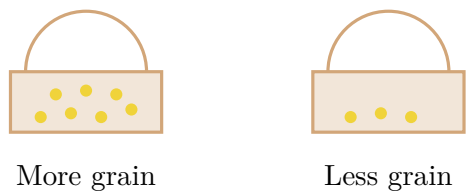
\begin{figure}[ht]
\centering

\begin{tikzpicture}[scale=1]

\draw[very thick,brown!70,fill=brown!20]
(-1,0) rectangle (1,0.8);

\draw[very thick,brown!70]
(-0.8,0.8) arc[start angle=180,end angle=0,radius=0.8cm];

\foreach \x/\y in {
-0.6/0.2,-0.2/0.25,0.2/0.2,0.6/0.3,
-0.4/0.5,0/0.55,0.4/0.5}
{
\fill[yellow!70!brown] (\x,\y) circle (0.08);
}

\draw[very thick,brown!70,fill=brown!20]
(3,0) rectangle (5,0.8);

\draw[very thick,brown!70]
(3.2,0.8) arc[start angle=180,end angle=0,radius=0.8cm];

\foreach \x/\y in {
3.5/0.2,3.9/0.25,4.3/0.2}
{
\fill[yellow!70!brown] (\x,\y) circle (0.08);
}

\node at (0,-0.6){More grain};

\node at (4,-0.6){Less grain};

\end{tikzpicture}

\caption{
Visual comparison reveals which basket contains more grain.
}

\label{fig:moregrain}

\end{figure}

The room became quiet. The teacher picked up the chalk and drew two baskets on the board.
\begin{quote}
\emph{"Imagine you are a farmer. Two baskets of grain are placed in front of you (Figure~\ref{fig:moregrain}). Which basket contains more grain?"}
\end{quote}
The answer came quickly.
\begin{quote}
\emph{"The basket on the left."}
\end{quote}
The teacher smiled.
\begin{quote}
\emph{"Exactly. When one collection is clearly larger, your eyes can tell."}
\end{quote}
He erased the drawing and replaced it with another.

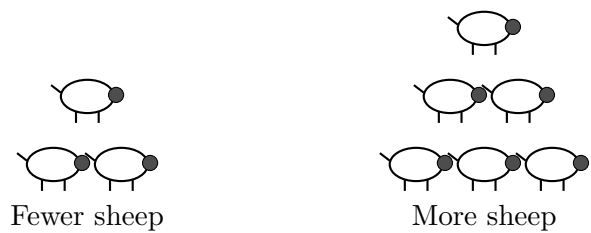
\begin{figure}[ht]
\centering

\begin{tikzpicture}[scale=1]

\newcommand{\sheep}[2]{

\draw[fill=white,thick]
(#1,#2) ellipse (0.35 and 0.22);

\draw[fill=black!70]
(#1+0.38,#2+0.02) circle (0.10);

\draw[thick]
(#1-0.15,#2-0.20)--(#1-0.15,#2-0.35);

\draw[thick]
(#1+0.15,#2-0.20)--(#1+0.15,#2-0.35);

\draw[thick]
(#1-0.35,#2+0.05)--(#1-0.48,#2+0.15);
}

\foreach \x/\y in {
0/0,0.9/0,
0.45/0.9}
{
\sheep{\x}{\y}
}

\foreach \x/\y in {
4.8/0,5.7/0,6.6/0,
5.25/0.9,6.15/0.9,
5.7/1.8}
{
\sheep{\x}{\y}
}

\node at (0.45,-0.7){Fewer sheep};

\node at (5.7,-0.7){More sheep};

\end{tikzpicture}

\caption{
A smaller and larger flock can be distinguished by sight.
}

\label{fig:sheepcomparison}

\end{figure}

\begin{figure}[ht]
\centering

\begin{tikzpicture}[scale=1]

\newcommand{\basket}[2]{

\draw[very thick,brown!70,fill=brown!20]
(#1-0.25,#2) rectangle (#1+0.25,#2+0.35);

\draw[very thick,brown!70]
(#1-0.18,#2+0.35)
arc[start angle=180,end angle=0,radius=0.18cm];

\fill[orange!80] (#1,#2+0.15) circle (0.07);
}

\foreach \x/\y in {
0/0,0.9/0,1.8/0,
0.45/1,1.35/1,
0.9/2,1.8/2}
{
\basket{\x}{\y}
}

\node at (0.9,-0.6){Collection A};

\basket{5.8}{2.4}
\basket{6.8}{2.0}
\basket{7.2}{1.0}
\basket{6.8}{0.2}
\basket{5.8}{-0.2}
\basket{4.8}{0.2}
\basket{4.4}{1.2}

\node at (5.8,-0.9){Collection B};

\end{tikzpicture}

\caption{
Different arrangements, same quantity.
}

\label{fig:samequantity}

\end{figure}
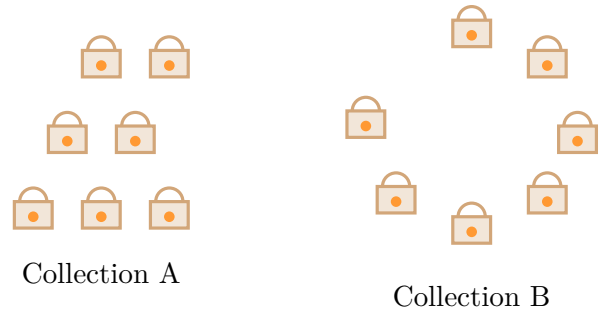

\begin{quote}
\emph{"Now imagine two shepherds returning with their flocks (Figure~\ref{fig:sheepcomparison}). Which flock has fewer sheep?"}
\end{quote}
Once again, the answer was immediate.
\begin{quote}
\emph{The flock on the left.}
\end{quote}
The teacher nodded.
\begin{quote}
\emph{"Good. When the difference is obvious, our eyes are enough."}
\end{quote}
Then he paused.
\begin{quote}
\emph{"But what happens when the difference is not obvious?"}
\end{quote}
A third drawing appeared on the board. The teacher pointed to the two collections (Figure~\ref{fig:samequantity}).


\begin{quote}
\emph{"Do these collections contain the same quantity?"}
\end{quote}
No one answered. The baskets were arranged differently. One spread across the board. The other formed a circle. Their appearance gave no obvious clue.
Hamid looked uncertain.
\begin{quote}
\emph{"They seem equal... but I cannot tell."}
\end{quote}
The teacher did not answer. Instead, he drew one more picture.


\begin{figure}[ht]
\centering

\begin{tikzpicture}[scale=1]


\foreach \x/\y in {
0/0,1.2/0,2.4/0,
0.6/1.4,1.8/1.4,
1.2/2.8,2.4/2.8}
{

\draw[very thick,brown!70,fill=brown!20]
(\x-0.35,\y) rectangle (\x+0.35,\y+0.45);

\draw[very thick,brown!70]
(\x-0.28,\y+0.45)
arc[start angle=180,end angle=0,radius=0.28cm];

\fill[orange!80] (\x-0.12,\y+0.18) circle (0.07);

\fill[red!60] (\x+0.12,\y+0.18) circle (0.07);

}

\node at (1.2,-0.7){Baskets};


\foreach \x/\y in {
5.5/0,6.7/0,7.9/0,
6.1/1.4,7.3/1.4,
6.7/2.8,7.9/2.8}
{

\fill[gray!65]
(\x,\y+0.15) ellipse (0.18 and 0.12);

}

\node at (6.7,-0.7){Pebbles};

\end{tikzpicture}

\caption{
One-to-one matching reveals equal quantities.
}

\label{fig:equalcollections}

\end{figure}
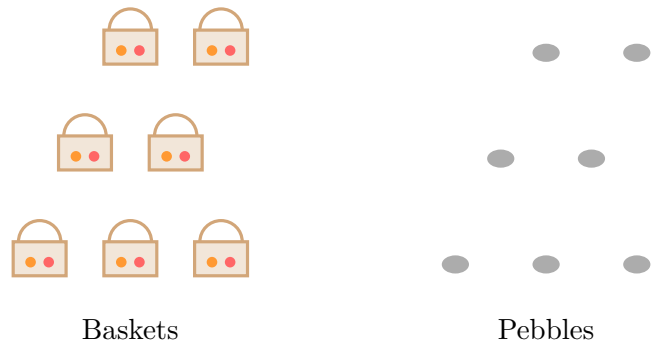

\begin{quote}
\emph{There was another way to compare them. People could match objects one by one. One basket with one pebble (Figure~\ref{fig:equalcollections}). One sheep with one stone.
If every object found a partner and nothing remained, the two collections had the same quantity. If something was left unmatched, one collection was larger.}
\end{quote}
Rajesh nodded.
\begin{quote}
\emph{"So they  only needed matching."}
\end{quote}
\begin{quote}
\emph{"Exactly,"} the teacher replied. \emph{"Today we call this \emph{one-to-one correspondence}.''}
\end{quote}
The classroom became quiet again. The idea was simple. But it answered a question people faced every day. There was still a problem. Matching worked only when the objects were present. Stones, shells, pebbles, and marks helped people compare quantities, but they could not easily preserve that information for later use. They could not tell someone, far away or years later, how many objects had once existed. 
\begin{quote}
\emph{How could people remember a quantity after the objects themselves were gone?}
\end{quote}


\section*{\textcolor{countinggreen}{The Birth of Counting: From Quantity to Number}}

The teacher looked around the classroom.
\begin{quote}
\emph{"Last time  we saw how people compared collections by matching objects one by one. It solved one problem. But another remained."}
\end{quote}


\begin{figure}[ht]
\centering

\begin{tikzpicture}[scale=1]

\newcommand{\sheep}[2]{
    \draw[fill=white,thick] (#1,#2) ellipse (0.30 and 0.18);
    \fill[black!70] (#1+0.33,#2+0.02) circle (0.08);
    \draw[thick] (#1-0.10,#2-0.18)--(#1-0.10,#2-0.33);
    \draw[thick] (#1+0.10,#2-0.18)--(#1+0.10,#2-0.33);
}

\foreach \x/\y in {
0/0,0.9/0,1.8/0,
0.45/0.8,1.35/0.8,
0.9/1.6,1.8/1.6}
{
\sheep{\x}{\y}
}

\node at (1.0,-0.7){Yesterday};

\node at (4.1,0.8){$\Longrightarrow$};

\node[font=\Huge] at (6.3,0.8){?};

\node at (6.3,-0.7){Today};

\end{tikzpicture}

\caption{ From physical objects to recorded quantity.}

\label{fig:memoryproblem}

\end{figure}
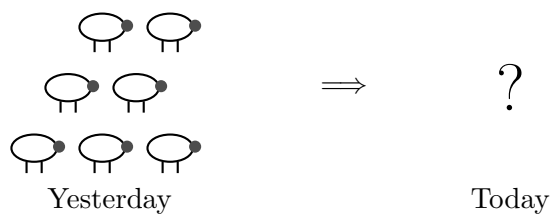

He drew a flock of sheep on the board (Figure~\ref{fig:memoryproblem}).

\begin{quote}
\emph{"Imagine a shepherd matching each sheep with a stone before sending the flock out to graze. The match confirms that none is missing. Later, the sheep wander away, leaving only the stones behind."}
\end{quote}

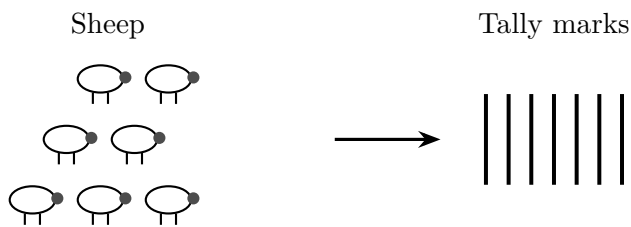
\begin{figure}[ht]
\centering

\begin{tikzpicture}[>=Stealth]

\newcommand{\sheep}[2]{
    \draw[fill=white,thick] (#1,#2) ellipse (0.30 and 0.18);
    \fill[black!70] (#1+0.33,#2+0.02) circle (0.08);
    \draw[thick] (#1-0.10,#2-0.18)--(#1-0.10,#2-0.33);
    \draw[thick] (#1+0.10,#2-0.18)--(#1+0.10,#2-0.33);
}

\foreach \x/\y in {
0/0,0.9/0,1.8/0,
0.45/0.8,1.35/0.8,
0.9/1.6,1.8/1.6}
{
\sheep{\x}{\y}
}

\node at (1.0,2.3){Sheep};

\draw[->,very thick] (4.0,0.8)--(5.4,0.8);

\foreach \x in {6.0,6.3,6.6,6.9,7.2,7.5,7.8}
{
\draw[line width=1.5pt] (\x,0.2)--(\x,1.4);
}

\node at (6.9,2.3){Tally marks};

\end{tikzpicture}

\caption{ One mark for each object preserves the quantity.}

\label{fig:tallymarks}

\end{figure}

One-to-one correspondence allowed comparison, but not memory. As communities grew and trade expanded, people needed a way to record quantities beyond the objects themselves. The answer was simple. People no longer preserved the objects; they preserved the quantity. One mark for each sheep (Figure~\ref{fig:tallymarks}) created a record that remained after the flock was gone.


\begin{figure}[ht]
\centering

\begin{tikzpicture}[scale=1]

\node[align=center] at (0,2.8){Seven\\Sheep};

\node[align=center] at (0,1.4){Seven\\Baskets};

\node[align=center] at (0,0){Seven\\Stones};

\draw[->,thick] (1.1,2.8)--(4.2,2.8);
\draw[->,thick] (1.1,1.4)--(4.2,1.4);
\draw[->,thick] (1.1,0)--(4.2,0);

\foreach \y in {2.8,1.4,0}
{
\foreach \x in {4.8,5.1,5.4,5.7,6.0,6.3,6.6}
{
\draw[line width=1.5pt] (\x,\y-0.35)--(\x,\y+0.35);
}
}

\end{tikzpicture}
\caption{ Tally marks represent quantity, not objects.}

\label{fig:abstraction}

\end{figure}
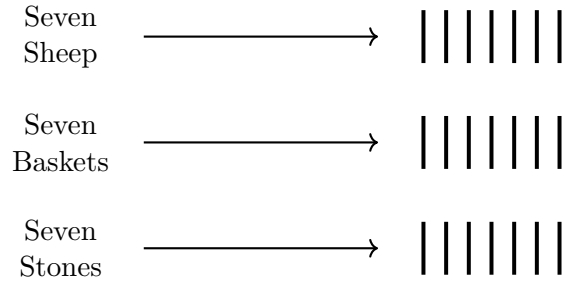

A deeper idea had appeared. The marks no longer belonged to sheep alone. They represented quantity itself (Figure~\ref{fig:abstraction}).


The teacher turned back to the class.
\begin{quote}
\emph{"Seven sheep, seven baskets of fruit, and seven stones are different collections. Yet all can produce the same tally record. "}
\end{quote}

\section*{\textcolor{arithmeticorange}{The First Arithmetic: Transforming Quantities}}

\begin{figure}[ht]
\centering

\begin{tikzpicture}[
    >=Stealth,
    node distance=1.8cm,
    every node/.style={font=\large}
]

\node[
draw,
rounded corners,
fill=blue!10,
minimum width=3.5cm,
minimum height=1cm
] (quantity)
{\textbf{Quantity}};


\node[
draw,
rounded corners,
fill=green!15,
above left=of quantity,
minimum width=3cm,
minimum height=1cm
] (add)
{\textbf{Combine}};

\node[
draw,
rounded corners,
fill=red!15,
below left=of quantity,
minimum width=3cm,
minimum height=1cm
] (sub)
{\textbf{Remove}
};

\node[
draw,
rounded corners,
fill=orange!15,
above right=of quantity,
minimum width=3cm,
minimum height=1cm
] (mul)
{\textbf{Combine equal  groups}
};

\node[
draw,
rounded corners,
fill=purple!15,
below right=of quantity,
minimum width=3cm,
minimum height=1cm
] (div)
{\textbf{Share equally}
};

\draw[->,thick] (quantity)--(add);
\draw[->,thick] (quantity)--(sub);
\draw[->,thick] (quantity)--(mul);
\draw[->,thick] (quantity)--(div);

\end{tikzpicture}

\caption{ Arithmetic arose from changing quantities.}

\label{fig:birthofarithmetic}

\end{figure}
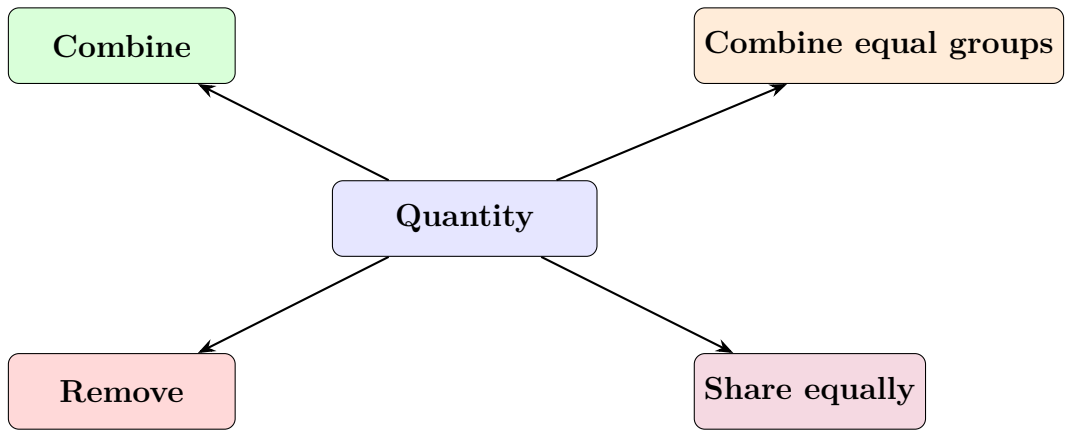

The teacher looked again at the tally marks.
\begin{quote}
\emph{"These marks represent quantity, not particular objects. But quantities rarely remain unchanged. They grow, shrink, combine, and divide. Describing those changes was the beginning of arithmetic (Figure~\ref{fig:birthofarithmetic})."}
\end{quote}


The tally marks were simple, but they recorded real actions: counting sheep, storing grain, trading goods, and dividing resources. People were already changing quantities. The marks only gave those changes a record.

\subsection*{Combining Quantities}

\begin{figure}[ht]
\centering

\begin{tikzpicture}[font=\large]

\node[align=center] at (0,1){Eastern field \\ Sheep};
\node at (0,0){||||};

\node at (2,0){combine};

\node[align=center] at (4,1){Western field\\ Sheep};
\node at (4,0){|||};

\node at (6,0){becomes};

\node[align=center] at (9,1){Total flock\\ Sheep};
\node at (9,0){|||||||};

\end{tikzpicture}

\caption{ \textbf{Addition} emerged from combining quantities.
}

\label{fig:additiontally}

\end{figure}
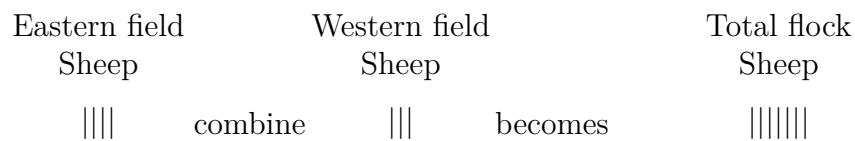

The teacher drew two groups of tally marks.
\begin{quote}
\emph{"Two shepherds record separate flocks. When the sheep return to one enclosure, what should happen to the records (Figure~\ref{fig:additiontally})?"}
\end{quote}
Hamid answered.
\begin{quote}
\emph{``Join the marks.''}
\end{quote}

\subsection*{Removing Quantities}

The teacher erased two tally marks.
\begin{quote}
\emph{"Some baskets of grain are used during the winter. What should happen to the record (Figure~\ref{fig:subtractiontally})?"}
\end{quote}
Fewer marks remained.

\begin{figure}[ht]
\centering

\begin{tikzpicture}[font=\large]

\node[align=center] at (0,1){baskets of grain};
\node at (0,0){|||||||};

\node at (2,0){remove};

\node[align=center] at (4,1){Used baskets};
\node at (4,0){||};

\node at (6,0){leaves};

\node[align=center] at (9,1){Remaining baskets};
\node at (9,0){|||||};

\end{tikzpicture}

\caption{ Removing quantities became the idea of \textbf{subtraction}.}

\label{fig:subtractiontally}

\end{figure}

\subsection*{Combining Equal Groups}
The teacher now drew three identical rows of tally marks.
\begin{quote}
\emph{"Three farmers harvest the same number of pumpkins (Figure~\ref{fig:multiplicationtally}). Must we count every mark, or can we use the pattern?"}
\end{quote}

\begin{figure}[ht]
\centering

\begin{tikzpicture}[font=\large]

\node[align=center] at (0,2){Farmer 1};
\node at (0,1.3){||||};

\node[align=center] at (0,0.5){Farmer 2};
\node at (0,-0.2){||||};

\node[align=center] at (0,-1){Farmer 3};
\node at (0,-1.7){||||};

\node at (3,0){combine equal groups};

\node[align=center] at (8,1){Total quantity};
\node at (8,0){||||||||||||};

\end{tikzpicture}

\caption{ Combining equal groups became the idea of \textbf{multiplication}.}

\label{fig:multiplicationtally}

\end{figure}

Repeated addition became a new way to view quantity.
\subsection*{Sharing Equally}

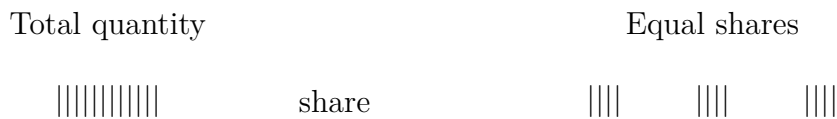
\begin{figure}[ht]
\centering

\begin{tikzpicture}[font=\large]

\node[align=center] at (0,1){Total quantity};
\node at (0,0){||||||||||||};

\node at (3,0){share};

\node[align=center] at (8,1){Equal shares};

\node at (8,0){|||| \qquad |||| \qquad ||||};

\end{tikzpicture}

\caption{
Sharing quantities equally became the idea of \textbf{division}.
}

\label{fig:divisiontally}

\end{figure}
The teacher drew one long row of tally marks.
\begin{quote}
\emph{"Now suppose the pumpkins are shared equally among three families. (Figure~\ref{fig:divisiontally}). What should happen to the record?"}
\end{quote}

The students separated the marks into three equal groups.

\section*{\textcolor{numeralplum}{From Numbers to Numerals: Giving Numbers a Written Form}}

\begin{figure}[ht]
\centering

\begin{tikzpicture}[node distance=2cm,>=Stealth]

\node[
draw,
rounded corners,
fill=blue!10,
minimum width=3.2cm,
minimum height=1cm
](objects){Five Objects};

\node[
draw,
rounded corners,
fill=green!10,
minimum width=3.2cm,
minimum height=1cm,
right=of objects
](number){Number: Five};

\node[
draw,
rounded corners,
fill=orange!15,
minimum width=3.8cm,
minimum height=1cm,
right=of number
](numerals){Numerals: $|||||,\;V,\;5$};

\draw[->,thick] (objects)--(number);
\draw[->,thick] (number)--(numerals);

\end{tikzpicture}

\caption{ Numbers represent quantity. Numerals represent numbers.
}

\label{fig:numbertonumeral}

\end{figure}
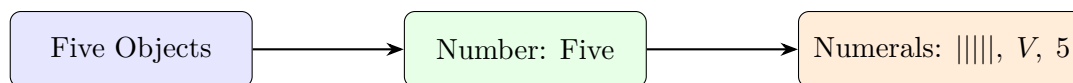

The teacher pointed again to the tally marks on the board.
\begin{quote}
\emph{"These marks let people record quantities, compare collections, and perform the earliest arithmetic. But imagine keeping track of one hundred sheep. How many marks would you need?"}
\end{quote}
\begin{quote}
\emph{"One hundred tally marks."}
\end{quote}
\begin{quote}
\emph{"Now imagine a merchant recording thousands of sacks of grain, or a ruler keeping accounts for an entire kingdom. Would tally marks still be practical?"}
\end{quote}
John answered first.
\begin{quote}
\emph{"It would take too long to write."}
\end{quote}
Rajesh added,
\begin{quote}
\emph{"And even longer to read."}
\end{quote}
\begin{quote}
\emph{"That follows. Tally marks work well for small quantities, but not for large ones. People needed a more efficient way to write numbers."}
\end{quote}
The teacher erased the board and wrote
\[
||||| \qquad\qquad \mathrm{V} \qquad\qquad 5
\]
The students  looked at  the three expressions.
\begin{quote}
\emph{"They don't look the same."}
\end{quote}
\begin{quote}
\emph{"The symbols are different, but the number is the same."}
\end{quote}
He pointed to each expression.
\begin{quote}
\emph{"Five tally marks record the quantity directly. The Romans wrote it as \(\mathrm{V}\). Today we write it as \(5\). Different numerals. One number."}
\end{quote}
Rajesh frowned.
\begin{quote}
\emph{"So are they different numbers?"}
\end{quote}
\begin{quote}
\emph{"No. A number is the abstract idea of quantity. A numeral is the symbol used to write it."}
\end{quote}
He drew the diagram shown in Figure~\ref{fig:numbertonumeral}.

The students looked again at
\[5\]
It no longer seemed like the number itself, only one way of writing it. Whether represented by tally marks, the Roman numeral \(\mathrm{V}\), or the modern digit \(5\), the number remained unchanged.
The teacher picked up the chalk again.
\begin{quote}
\emph{"One problem still remained. If every number needed its own symbol, the system would soon become unmanageable."}
\end{quote}
He wrote \[ 1,\;2,\;3,\;\ldots,\;9 \] on the board.
\begin{quote}
\emph{"These are not just symbols to memorize. They are the building blocks of mathematical language, much like letters of the alphabet. With just a small set of numerals, people could represent quantities far beyond nine."}
\end{quote}

\section*{\textcolor{indianmathgold}{The Indian Breakthrough: Place Value and the Birth of Zero}}

\begin{figure}[ht]
\centering

\begin{tikzpicture}[
    every node/.style={font=\large}
]


\node at (1.8,4) {\textbf{Ten sticks}};

\foreach \x in {0,0.35,0.7,1.05,1.4,1.75,2.1,2.45,2.8,3.15}
{
    \draw[very thick,brown!80]
    (\x,2.8)--(\x,3.8);
}

\draw[->,very thick]
(3.8,3.3)--(5,3.3);

\node at (4.4,3.8) {tie together};

\node at (6.5,4) {\textbf{One bundle}};

\draw[very thick,red!70]
(6.5,3.3) circle (0.3);

\node at (6.5,3.3) {};


\node at (2.5,1.5) {\textbf{Ten bundles}};

\foreach \x in {0,0.65,1.3,1.95,2.6,3.25,3.9,4.55,5.2,5.85}
{
    \draw[very thick,red!70]
    (\x,0.4) circle (0.3);
}

\draw[->,very thick]
(6.7,0.4)--(8,0.4);

\node at (7.35,0.9) {combine};

\node at (9.5,1.5) {\textbf{One larger unit}};

\draw[very thick,blue!70,rounded corners]
(8.8,-0.2) rectangle (10.2,1);

\node at (9.5,0.4) {};

\end{tikzpicture}

\caption{
Grouping creates larger units from repeated quantities.
}
\label{fig:bundlesoftens}
\end{figure}
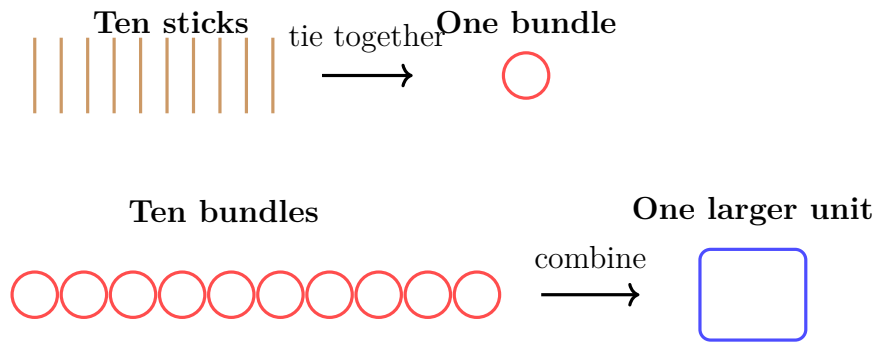

The teacher erased the board.
\begin{quote}
\emph{"Imagine you are a merchant with three hundred and twenty-three sacks of grain. Using only these nine symbols, how would you write that number?"}
\end{quote}
The room became quiet.
\begin{quote}
\emph{"Could we just repeat the symbols?"}
\end{quote}
The teacher smiled.
\begin{quote}
\emph{"That works for small numbers. But for thousands of sacks, the records would become too long and difficult to manage. People needed a shorter way to write large quantities."}
\end{quote}
He then picked up a handful of sticks.
\begin{quote}
\emph{"Ten sticks can be grouped into one larger unit—a bundle of ten. We no longer need to count each stick; we can count bundles instead. Ten such bundles can then be grouped into an even larger unit—a hundred (Figure~\ref{fig:bundlesoftens})."}
\end{quote}


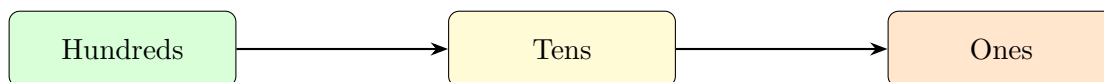
\begin{figure}[ht]
\centering

\begin{tikzpicture}[>=Stealth,node distance=2.8cm]

\node[draw,rounded corners,fill=green!15,
minimum width=3cm,minimum height=1cm] (H)
{Hundreds};

\node[draw,rounded corners,fill=yellow!20,
minimum width=3cm,minimum height=1cm,
right=of H] (T)
{Tens};

\node[draw,rounded corners,fill=orange!20,
minimum width=3cm,minimum height=1cm,
right=of T] (O)
{Ones};

\draw[->,thick] (H)--(T);
\draw[->,thick] (T)--(O);

\end{tikzpicture}

\caption{
Place value organizes the quantities into hundreds, tens, and ones.
}
\label{fig:placevalue}
\end{figure}

The students immediately saw the advantage. Large quantities no longer required endless marks. A few symbols could now represent many objects.
Then  the teacher wrote:
\begin{quote}
\emph{"For three hundred and twenty-three sacks of grain, how many groups of one hundred?"}
\end{quote}
\begin{quote}
\emph{``Three.''}
\end{quote}
\begin{quote}
\emph{``Groups of ten?''}
\end{quote}
\begin{quote}
\emph{``Two.''}
\end{quote}
\begin{quote}
\emph{``And single sacks?''}
\end{quote}
\begin{quote}
\emph{``Three.''}
\end{quote}
\begin{quote}
\emph{``So here quantity can be described by hundreds, tens, and ones.''}
\end{quote}
He drew three boxes (Figure~\ref{fig:placevalue}).
Then he wrote
\[323\]
\begin{quote}
\emph{"The first digit means three hundreds, the second means two tens, and the last means three ones."}
\end{quote}
He pointed to the two 3s.
\begin{quote}
\emph{"They look the same. Why are they different?"}
\end{quote}
\begin{quote}
\emph{"Because they are in different places,"}
\end{quote}
Rajesh replied.
\begin{quote}
\emph{"Exactly. A digit's value comes from its position. This idea, called \emph{place value}, became one of the most important advances in mathematics and was developed in India."}
\end{quote}
The teacher erased the middle digit and wrote
\[ 3 \qquad\qquad 3 \]
\begin{quote}
\emph{"Now suppose there are three hundreds, no tens, and three ones. How should we write it?"}
\end{quote}
John frowned.
\begin{quote}
\emph{"Writing 3 nothing 3 would mean three hundreds, no tens and three ones."}
\end{quote}
\begin{quote}
\emph{"Right. We need a symbol to hold the empty place."}
\end{quote}
He wrote
\[303\]
\begin{quote}
\emph{"Zero keeps the place when a value is missing. Without it, place value cannot work."}
\end{quote}
Rajesh raised his hand.
\begin{quote}
\emph{"Was zero created only for an empty place?"}
\end{quote}
\begin{quote}
\emph{"That was its first use."}
\end{quote}
The teacher wiped the board.
\begin{quote}
\emph{"If every sheep leaves a field, how many remain?"}
\end{quote}
\begin{quote}
\emph{"None."}
\end{quote}
He wrote
\[0\]
\begin{quote}
\emph{"Now even nothing has a symbol. Zero is not only a placeholder; it also represents the absence of quantity."}
\end{quote}
The teacher stepped back.
\begin{quote}
\emph{"With place value and zero, people could write any number. But another question remained: how do we represent what we do not have—what we owe, lose, or lack?"}
\end{quote}
The students looked at one another. The journey beyond counting had begun.

\section*{\textcolor{crimson}{The Merchant's Problem: The Birth of Negative Numbers}}

The teacher drew a merchant's record.
\begin{center}
\begin{tabular}{|l|}
\hline
Merchant A\\
\hline
Owns three gold coins\\
\hline
\end{tabular}
\end{center}
\begin{quote}
\emph{"Can our numbers describe this merchant?"}
\end{quote}
\begin{quote}
\emph{"Three."}
\end{quote}
The teacher changed the record.
\begin{center}
\begin{tabular}{|l|}
\hline
Merchant B\\
\hline
Owes three gold coins\\
\hline
\end{tabular}
\end{center}

The room became quiet.
\begin{quote}
\emph{"What number describes this merchant?"}
\end{quote}
John hesitated.
\begin{quote}
\emph{"Not three."}
\end{quote}
The teacher nodded.
\begin{quote}
\emph{"Yes. Existing numbers described what people had, but not what they owed."}
\end{quote}
He wrote:
\[
\text{Wealth}\qquad\qquad\text{Debt}
\]
\begin{quote}
\emph{"These are opposites. The same amount can have opposite meanings: profit and loss, forward and backward, wealth and debt. Mathematics needed numbers for both directions."}
\end{quote}
John raised his hand.
\begin{quote}
\emph{"Then what is the opposite of debt?"}
\end{quote}

\begin{quote}
\emph{"Wealth. The opposite of an opposite returns us to the original quantity."}
\end{quote}
He wrote:
\[
\text{opposite of opposite quantity}\; equals \; \text{original quantity}
\]
\begin{quote}
\emph{"This idea became the foundation of signed numbers."}
\end{quote}
The teacher stepped back.
\begin{quote}
\emph{"Counting was no longer enough. Mathematics needed numbers that could represent opposite quantities."}
\end{quote}
A new question appeared:
\begin{quote}
\emph{"How should these numbers be written, and how should arithmetic work with them?"}
\end{quote}
The journey into negative numbers had begun.

\section*{\textcolor{symbolicindigo}{From Commerce to Symbolic Mathematics}}

The teacher returned to the board.
\begin{quote}
\emph{"We have seen why merchants needed numbers that could represent both wealth and debt. But a new question appeared: how should these quantities be written?"}
\end{quote}
John raised his hand.
\begin{quote}
\emph{"A merchant could write `owns three coins' or `owes three coins."}
\end{quote}
The teacher nodded.
\begin{quote}
\emph{"That works for one record. But what about thousands of transactions every day? Would words be the best way to describe every change?"}
\end{quote}
The students shook their heads.
\begin{quote}
\emph{"It would take too long."}
\end{quote}
Rajesh added,
\begin{quote}
\emph{"And people might describe the same situation in different ways."}
\end{quote}
\begin{quote}
\emph{"Good observation,"}
\end{quote}
said the teacher.
\begin{quote}
\emph{"Mathematics was growing, but ordinary language could not keep up. People needed a shorter and clearer way to express ideas."}
\end{quote}
He picked up the chalk.
\begin{quote}
\emph{"Over time, that language began to take shape."}
\end{quote}
He wrote:
\[=\]
\begin{quote}
\emph{"This symbol represented equality."}
\end{quote}
Then he added:
\[
+\qquad -\qquad \times\qquad \div
\]
\begin{quote}
\emph{"These symbols represented actions people had performed for centuries: joining, removing, grouping, and sharing quantities."}
\end{quote}
Finally, he wrote:
\[+2\qquad\qquad-3\]
The students looked at the expressions.
\begin{quote}
\emph{"The symbols are small,"}
\end{quote}
John said,
\begin{quote}
\emph{"but they immediately show whether a quantity represents a gain or a loss."}
\end{quote}
The teacher smiled.
\begin{quote}
\emph{"That is the point. A single symbol could now express an idea that once required a sentence."}
\end{quote}
He turned to the class.
\begin{quote}
\emph{"But opposite situations are not limited to trade. Where else do we see them?"}
\end{quote}
The students answered:
\begin{quote}
\emph{"Moving forward and backward."}
\end{quote}
\begin{quote}
\emph{"Going upward and downward."}
\end{quote}
The teacher nodded.
\begin{quote}
\emph{"An idea that began with debts and accounts became a way to describe direction and opposite states as well."}
\end{quote}
He continued,
\begin{quote}
\emph{"Astronomers used similar ideas to describe positions and movement. Scientists used them to represent increases and decreases. What began as a need in commerce became part of the language of mathematics."}
\end{quote}
The teacher stepped away from the board. The board was filled with familiar symbols:
\[
\qquad
+
\qquad
-
\qquad
\times
\qquad
\div
\]
They looked simple. They were not. Each symbol carried the memory of problems people once faced and the solutions they created. Each one marked a step in humanity’s effort to understand quantity and change.

\section*{\textcolor{tealgreen}{Extending Arithmetic: When Operations Acquired New Meanings}}

The teacher looked around the classroom.
\begin{quote}
\emph{"We have travelled a long way. We began by comparing baskets and flocks. Then came counting, numerals, place value, and zero. Later, numbers grew to represent opposites---wealth and debt, gain and loss, forward and backward."}
\end{quote}
He paused before continuing.
\begin{quote}
\emph{"Now another question arises. Once negative numbers entered mathematics, could arithmetic still mean exactly what it had meant before?"}
\end{quote}
John raised his hand.
\begin{quote}
\emph{"Why not? Addition is still addition, and multiplication is still multiplication."}
\end{quote}
The teacher smiled.
\begin{quote}
\emph{"It seems that way. But arithmetic began with physical actions---adding sheep, removing grain, and sharing supplies. Negative numbers changed what a quantity could represent. Numbers no longer described only objects. They could also describe opposite effects."}
\end{quote}
The room fell silent. The teacher drew two records on the board.
\begin{center}
\begin{tabular}{|c|}
\hline
Merchant A\\
\hline
Owns Rs 800\\
\hline
\end{tabular}
\hspace{0.5cm}
\begin{tabular}{|c|}
\hline
Merchant B\\
\hline
Owns Rs 500\\
\hline
\end{tabular}
\end{center}

\vspace{0.4cm}

\begin{center}
\begin{tabular}{|c|}
\hline
Merchant A\\
\hline
Owns Rs 800\\
\hline
\end{tabular}
\hspace{0.5cm}
\begin{tabular}{|c|}
\hline
Merchant B\\
\hline
Owes Rs 500\\
\hline
\end{tabular}
\end{center}

The teacher pointed to the second pair of records.
\begin{quote}
\emph{"Now what happens?"}
\end{quote}
Rajesh studied the records.
\begin{quote}
\emph{"The debt reduces the wealth."}
\end{quote}
The teacher nodded.
\begin{quote}
\emph{"Exactly. Addition is no longer just about combining quantities. It combines their effects. Wealth pushes the balance in one direction. Debt pushes it in the opposite direction."}
\end{quote}
John nodded.
\begin{quote}
\emph{"So gains reinforce gains. Debts reinforce debts. But gains and debts cancel each other."}
\end{quote}
\begin{quote}
\emph{"Well said. The operation has not changed. Only the meaning of the quantities has expanded."}
\end{quote}
The teacher erased the board and wrote:

\begin{center}
\begin{tabular}{|c|}
\hline
Merchant\\
\hline
Owes Rs 500\\
\hline
\end{tabular}
\end{center}

\begin{quote}
\emph{"Suppose a friend pays this debt. What happened? Did the merchant receive Rs 500?"}
\end{quote}

The students shook their heads.

\begin{quote}
\emph{"No."}
\end{quote}

John answered.

\begin{quote}
\emph{"Then why is the merchant in a better position?"}
\end{quote}

Rajesh replied,

\begin{quote}
\emph{"Because the debt is gone."}
\end{quote}

The teacher nodded.

\begin{quote}
\emph{"That is correct. Nothing was added to the merchant's hands, yet something negative disappeared."}
\end{quote}

He turned back to the class.

\begin{quote}
\emph{"What has the same effect as removing a debt?"}
\end{quote}

The room became quiet. After a moment, John answered.

\begin{quote}
\emph{"Adding wealth."}
\end{quote}

The teacher smiled.

\begin{quote}
\emph{"Correct. Removing a debt has the same effect as adding an equal amount of wealth. Likewise, removing wealth has the same effect as adding debt."}
\end{quote}

He wrote

\[
\text{opposite}
\]

on the board.

\begin{quote}
\emph{"This changed the way subtraction was understood. Subtracting a quantity could now be viewed as adding its opposite. Addition and subtraction were no longer separate ideas. They became two closely connected ways of describing the same change."}
\end{quote}

The students looked back at the board. Arithmetic was not becoming more complicated. It was becoming more expressive.

The teacher cleared the board and wrote

\[
2\times(\text{wealth of Rs 3})
\]

He turned toward the class.

\begin{quote}
\emph{"What does this expression mean?"}
\end{quote}

The students answered together.

\begin{quote}
\emph{"It means two amounts of wealth, each worth Rs 3, so the total wealth is Rs 6."}
\end{quote}

The teacher nodded.

\begin{quote}
\emph{"Indeed. Multiplication describes repeated groups. The first number tells us how many times the quantity is taken."}
\end{quote}

He paused.

\begin{quote}
\emph{"Now suppose the quantity points in the opposite direction."}
\end{quote}

He wrote

\[
2\times(\text{debt of Rs 3})
\]

on the board.

\begin{quote}
\emph{"Can we still understand this as repeated groups?"}
\end{quote}

Rajesh answered,

\begin{quote}
\emph{"Yes. It means two debts of Rs 3, so the total debt is Rs 6."}
\end{quote}

The teacher smiled.

\begin{quote}
\emph{"Precisely. Repetition still works. The quantity is simply pointing in the negative direction."}
\end{quote}

Then he wrote

\[
(-3)\times(\text{2 steps to the east})
\]

The room became quiet. John looked at the board.

\begin{quote}
\emph{"Repeated groups no longer seem to explain it."}
\end{quote}

The teacher nodded.

\begin{quote}
\emph{"This is the key point. The old interpretation has reached its limit. Before negative numbers, multiplication answered only one question: \emph{How many times?} Once numbers could point in opposite directions, multiplication had to answer another question: \emph{What happens to direction?}"}
\end{quote}

He drew a horizontal line across the board.

\begin{center}
Negative direction (West)\hspace{3cm} Positive direction (East)
\end{center}


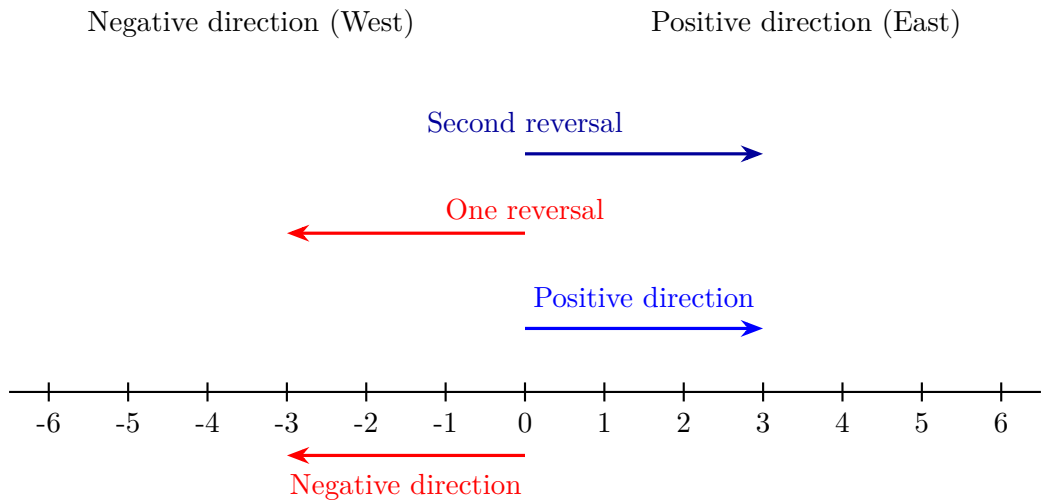
\begin{figure}[ht]
\centering

\begin{tikzpicture}[>=Stealth,scale=1.05]

\draw[-,thick] (-6.5,0)--(6.5,0);

\foreach \x in {-6,-5,-4,-3,-2,-1,0,1,2,3,4,5,6}
{
    \draw[thick] (\x,0.12)--(\x,-0.12);
    \node[below] at (\x,-0.15){\x};
}

\draw[->,very thick,blue] (0,0.8)--(3,0.8);
\node[blue] at (1.5,1.2){Positive direction};

\draw[->,very thick,red] (0,-0.8)--(-3,-0.8);
\node[red] at (-1.5,-1.2){Negative direction};

\draw[->,very thick,red] (0,2)--(-3,2);
\node[red] at (0,2.3){One reversal};

\draw[->,very thick,blue!60!black] (0,3)--(3,3);
\node[blue!60!black] at (0,3.4){Second reversal};

\end{tikzpicture}

\caption{A negative sign reverses direction; two reversals restore the original direction.}

\label{fig:directionreversal}

\end{figure}

The teacher pointed to the arrows in Figure~\ref{fig:directionreversal}.

\begin{quote}
\emph{"Look carefully. A negative sign reverses direction. One reversal changes the direction once. A second reversal changes it again, bringing us back to the original direction."}
\end{quote}

He turned back to the class.

\begin{quote}
\emph{"A signed quantity has two parts: magnitude and direction. The magnitude tells us how much. The sign tells us which way."}
\end{quote}

John looked at the diagram.

\begin{quote}
\emph{"So the sign tells us the direction of the quantity."}
\end{quote}

\begin{quote}
\emph{"Right,"}
\end{quote}

said the teacher.  He wrote

\[
(+2)\times(\text{3 steps to the east})
\]

\begin{quote}
\emph{"A positive multiplier does not change the direction. It simply repeats the quantity. Since the quantity already points west, repeating it twice gives six steps to the west."}
\end{quote}

Then he wrote

\[
(-3)\times(\text{2 steps to the east})
\]

\begin{quote}
\emph{"Now the multiplier is negative. What should it do?"}
\end{quote}

The class remained silent. After a moment, Rajesh answered,

\begin{quote}
\emph{"It should reverse the direction."}
\end{quote}

The teacher smiled.

\begin{quote}
\emph{"Yes, that is the idea. A negative multiplier acts as a reversal. It changes two steps east into two steps west. The reversed quantity is then repeated three times, giving six steps west."}
\end{quote}

He wrote
\[
\begin{aligned}
(-3)\times(\text{2 steps to the east})
&=(\text{2 steps to the west}) \\
&\quad+(\text{2 steps to the west}) \\
&\quad+(\text{2 steps to the west}) \\
&=\text{6 steps to the west}.
\end{aligned}
\]

The students nodded.

The new interpretation remained consistent with everything they had learned about multiplication. Only one new idea had been added: a negative multiplier reverses direction.

The teacher wrote

\[
(+2)\times(+3)
\]

on the board.

Rajesh answered immediately,

\begin{quote}
\emph{"The answer is \(+6\)."}
\end{quote}

The teacher asked,

\begin{quote}
\emph{"Why?"}
\end{quote}

Rajesh replied,

\begin{quote}
\emph{"The multiplier \(+2\) leaves the direction unchanged. The quantity \(+3\) stays positive and is repeated twice:
\[
(+3)+(+3)=+6.
\]"}
\end{quote}

The teacher nodded and wrote

\[
(+2)\times(-3).
\]

John answered,

\begin{quote}
\emph{"The answer is \(-6\)."}
\end{quote}

\begin{quote}
\emph{"Good observation. The multiplier is still positive, so the direction does not change. The quantity already points in the negative direction, so repeating it twice gives
\[
(-3)+(-3)=-6.
\]"}
\end{quote}

The teacher wrote one more example:

\[
(-3)\times(+2).
\]

Hamid answered,

\begin{quote}
\emph{"The answer is \(-6\)."}
\end{quote}

The teacher smiled.

\begin{quote}
\emph{"Can you explain why?"}
\end{quote}

Hamid replied,

\begin{quote}
\emph{"The multiplier is negative, so it reverses the direction of the quantity. The quantity \(+2\) becomes \(-2\), and repeating it three times gives
\[
(-2)+(-2)+(-2)=-6.
\]"}
\end{quote}

The teacher nodded.

\begin{quote}
\emph{"You got it. A positive multiplier preserves direction. A negative multiplier reverses it."}
\end{quote}

Finally, he wrote

\[
(-2)\times(-3).
\]

The class looked at the expression again. This time, it no longer seemed mysterious. The teacher pointed to the quantity \((-3)\).

\begin{quote}
\emph{"This quantity already points in the negative direction. What happens when a negative multiplier acts on it?"}
\end{quote}

After a brief pause, Rajesh answered,

\begin{quote}
\emph{"The multiplier reverses its direction. So \(-3\) becomes \(+3\)."}
\end{quote}

The teacher nodded.

\begin{quote}
\emph{"Correct. Once the direction is reversed, the positive quantity \(+3\) is repeated twice:
\[
(+3)+(+3)=+6.
\]"}
\end{quote}

He completed the calculation:

\[
(-2)\times(-3)=+6.
\]

Then he slowly wrote

\[
\boxed{(-)\times(-)=+}
\]

The classroom remained silent for a moment. Finally, the teacher said,
\begin{quote}
\emph{"This result was not invented to create another rule. It emerged naturally as the number system expanded. Once numbers could represent opposite directions, arithmetic had to remain consistent with that idea."}
\end{quote}

He looked back at the path they had followed.

\begin{quote}
\emph{"We began by comparing collections. We learned to count. We created numerals, place value, and zero. Then numbers grew to represent opposites. Finally, arithmetic itself expanded so that these new numbers could behave consistently."}
\end{quote}

\section*{\textcolor{rose}{Conclusion: From Counting Objects to Understanding Ideas}}

The expression

\[
(-2)\times(-3)=+6
\]

looks simple. But behind it lies a long history of human attempts to understand quantity, change, and direction. The journey began with comparisons. People compared baskets of fruit, flocks of sheep, and collections of objects. These comparisons created the need for counting. Counting led to numbers. Numbers required symbols. Larger societies needed better ways to record quantities, leading to place value. Place value required a symbol for an empty position, and zero transformed the number system. With zero came the possibility of describing quantities below zero. Negative numbers entered mathematics. Then arithmetic itself had to expand. Addition became more than combining objects; it became a way to combine effects. Subtraction became connected with opposites. Multiplication had to account not only for size but also for direction. The rule
[
\boxed{(-)\times(-)=+}
]
was not a shortcut to memorize. It was a consequence. The idea of direction explains why two reversals return to the original direction. The laws of arithmetic confirm that this interpretation keeps the number system consistent. 

\section*{Acknowledgement}

The author gratefully acknowledges Prof.\ Bharat Madhusudan Deshpande for his insightful discussions, constructive suggestions, and valuable comments, which greatly improved the quality and presentation of this work.

\end{document}